\documentclass[11pt,a4paper]{article}%
\usepackage{amscd}
\usepackage{amsmath}
\usepackage{graphicx}
\usepackage{amsfonts}
\usepackage{amssymb}%
\newtheorem{theorem}{Theorem}[section]
\newtheorem{lemma}[theorem]{Lemma}

\newtheorem{remark}[theorem]{Remark}

\numberwithin{equation}{section}
\def\qed { \vskip 0pt \hfill \hbox{\vrule height 5pt width 5pt depth 0pt} \vskip 10pt}

\begin{document}

\title{Self-similar solutions for active scalar equations in Fourier-Besov-Morrey
spaces }
\author{\textbf{Lucas C. F. Ferreira}\\{\small Universidade Estadual de Campinas, IMECC- Departamento de
Matem\'atica,} \\{\small Rua S\'ergio Buarque de Holanda, 651, CEP 13083-859, Campinas-SP,
Brazil.} \\{\small \texttt{email:\ lcff@ime.unicamp.br}} \vspace{1cm}\\\textbf{Lidiane S. M. Lima}\\{\small Universidade Estadual de Campinas, IMECC- Departamento de
Matem\'atica,} \\{\small Rua S\'{e}rgio Buarque de Holanda, 651, CEP 13083-859, Campinas-SP,
Brazil.} \\{\small \texttt{email:\ lidynet@hotmail.com}}}
\date{}
\maketitle

\begin{abstract}
We are concerned with a family of dissipative active scalar equation with
velocity fields coupled via multiplier operators that can be of high-order. We
consider sub-critical values for the fractional diffusion and prove global
well-posedness of solutions with small initial data belonging to a framework
based on Fourier transform, namely Fourier-Besov-Morrey spaces. Since the
smallness condition is with respect to the weak norm of this space, some
initial data with large $L^{2}$-norm can be considered. Self-similar solutions
are obtained depending on the homogeneity of the initial data and couplings.
Also, we show that solutions are asymptotically self-similar at infinity. Our
results can be applied in a unified way for a number of active scalar PDEs
like 1D models on dislocation dynamics in crystals, Burguer's equations, 2D
vorticity equation, 2D generalized SQG, 3D magneto-geostrophic equations,
among others.



\end{abstract}

%%%%%%%%%%%%%%%%%%%%%%%%%%%%%%%%%%%%%%%%%%%%%%%%%%%%%%%%%%%%%%%%%%%%%%%%%%%%%%%%%%%%%%%%%%%%%%%%%%%%%%%%%%%%%%%%%%%%

\pagestyle{myheadings} \markright{Dissipative active scalar equations}

\section{Introduction}

\baselineskip17pt\bigskip

\hspace{0.5cm} Frameworks based on the Fourier transform give an elegant way
of handling PDEs and may reveal underlying mathematical and physical aspects
connected with low and high frequencies. Our main intent is to obtain a theory
of self-similar solutions and stability in such type of framework that can be
applied in a unified way for a number of active scalar PDEs.

For that matter, we consider a family of dissipative active scalar equations
with velocity fields coupled via certain multiplier operators, which includes
PDEs with hamiltonian or gradient flow structure. Examples of them are given
below together with the corresponding literature. In fact, nice physical
models arise in dimension $n=1,2,3$. That family was introduced in
\cite{Chae2} (see also \cite{Chae1}) and reads as
\begin{equation}%
\begin{cases}
\frac{\partial\theta}{\partial t}+\kappa\left(  -\Delta\right)  ^{\gamma
}\theta+\nabla_{x}\cdot(u\theta)=0, & \qquad x\in\mathbb{R}^{n}%
\ ,\ t>0,\ \\[3mm]%
\theta(x,0)=\theta_{0}(x), & \qquad x\in\mathbb{R}^{n},\
\end{cases}
\label{dase}%
\end{equation}
where $n\geq1$, $\gamma>1/2,$ and the fractional dissipation is defined by
$\ \widehat{(-\Delta)^{\gamma}f}](\xi)=|\xi|^{2\gamma}\hat{f}(\xi).$ Since we
are concerned with the dissipative case, we assume $\kappa=1$ for the sake of simplicity.

The velocity field $u$ is coupled to the scalar $\theta$ via the multiplier
vector operator
\begin{equation}
u=P[\theta]=(u_{1},u_{2},...,u_{n})\label{velocity}%
\end{equation}
where
\begin{equation}
u_{k}=\sum_{j=1}^{n}a_{jk}\mathcal{R}_{j}\Lambda^{-1}P_{j}[\theta
],\;\;\text{for}\ 1\leq k\leq n,\label{uks}%
\end{equation}
\bigskip$\Lambda=(-\Delta)^{\frac{1}{2}}$, $\mathcal{R}_{j}=-\partial
_{j}(-\Delta)^{-\frac{1}{2}}$ is the $j$-th Riesz transform, $a_{jk}$'s are
constant and
\begin{equation}
\widehat{P_{j}[\theta]}(\xi)=P_{j}(\xi)\widehat{\theta}(\xi)\text{.}%
\label{Pjotas}%
\end{equation}
In the case $n=1,$ $\mathcal{R}_{j}$'s should be understood as the Hilbert
transform
\[
\mathcal{H}(u)(x)=\frac{1}{\pi}P.V.\int_{-\infty}^{\infty}\frac{u(y)}{x-y}dy.
\]
Let $\mathcal{A}$ be the set of all $f\in\mathcal{S}^{\prime}(\mathbb{R}^{n})$
such that $\widehat{f}|_{U}\in\mathcal{D}^{\prime}(U)$ is a complex Radon
measure on $U$, for all bounded open $U$, and the total variation $|\hat{f}|$
(extended to $\mathbb{R}^{n}$) is a tempered measure, that is, $\int
_{\mathbb{R}^{n}}(1+\left\vert x\right\vert )^{-N}d|\hat{f}|<\infty$, for some
$N\geq0$. We assume that $P_{j}$'s are measurable functions such that
\begin{equation}
|P_{j}(\xi)|\leq C|\xi|^{\beta}\text{, a.e. }\xi\in\mathbb{R}^{n}%
\text{,}\label{Pj-cond}%
\end{equation}
for all $1\leq j\leq n$, where the parameter $\beta\in\lbrack0,n+1)$ with
$\beta<2\gamma$; and so $P_{j}[\cdot]$ makes sense from $\mathcal{A}$ to
$\mathcal{S}^{\prime}(\mathbb{R}^{n})$. In view of (\ref{uks}) and
(\ref{Pjotas}), we can write the Fourier transform of the velocity field $u$
as
\begin{equation}
\hat{u}(\xi)=P(\xi)\hat{\theta}(\xi)\label{velocityfourier}%
\end{equation}
where $P(\xi)=$ $[\tilde{P}_{1},..,\tilde{P}_{n}]$ with
\[
\tilde{P}_{k}(\xi)=\sum_{j=1}^{n}a_{jk}\frac{i\xi_{j}}{\left\vert
\xi\right\vert ^{2}}P_{j}(\xi),\text{ }1\leq k\leq n,
\]
and so
\begin{equation}
|\hat{u}(\xi)|\leq C|\xi|^{\beta-1}|\hat{\theta}(\xi)|.\label{est velocity}%
\end{equation}

Assuming that $P_{j}$'s are homogeneous of degree $\beta,$ (\ref{dase}) has
the scaling map
\begin{equation}
\theta\rightarrow\theta_{\lambda}=\lambda^{2\gamma-\beta}\theta(\lambda
x,\lambda^{2\gamma}t),\text{ for all }\lambda>0, \label{scaling}%
\end{equation}
which induces the scaling for the initial condition
\begin{equation}
\theta_{0}\rightarrow\lambda^{2\gamma-\beta}\theta_{0}(\lambda x).
\label{scaling initial}%
\end{equation}
Even when $P_{j}$'s are not homogeneous, the map (\ref{scaling}) works well as
an intrinsic scaling for (\ref{dase})-(\ref{velocity}) in the sense that it is
useful to identify suitable functional spaces for global existence and
provides a notion of criticality. Indeed, we have three basic cases:
sub-critical $\beta<2\gamma,$ critical $\beta=2\gamma,$ and super-critical
$\beta>2\gamma$, which correspond to the sign of the exponent in
(\ref{scaling initial}).

From (\ref{est velocity}), one can see that (\ref{dase})-(\ref{velocity})
takes into account the effect of couplings with positive order when $\beta>1$,
which are named here as high-order ones. They behave morally like a positive
derivative of $(\beta-1)$-order, and produce more difficulties in comparison
with SQG ($\beta=1$) and $\beta<1,$ at least in the context of $H^{s}$ and
$L^{p}$-theory.

Here we consider a functional setting based on Fourier transform in which the
existence theory works in a unified way for both low ($\beta<1$) and
high-order couplings. More precisely, we prove global well-posedness and
stability results for (\ref{dase})-(\ref{velocity}) with small initial data
belonging to the framework of homogeneous Fourier-Besov-Morrey spaces
$\mathcal{FN}_{p,\mu,\infty}^{s}$ (see (\ref{FB-norm1}) for the definition),
where $s=n-\frac{n-\mu}{p}-(2\gamma-\beta)$, $\frac{n-\mu}{n+\beta+1-4\gamma
}<p\leq\infty$, $0\leq\mu<n$, $\gamma\neq1/2$ and $0\leq\beta<2\gamma
<\frac{n+\beta+1}{2}$ (see Theorem \ref{teoglobal}). In view of the norm
(\ref{FB-norm1}), our results allow to take some large data in $H^{s}$ and
$L^{p}$-spaces (see Remark \ref{Rem2}), and homogeneous functions of degree
$-(2\gamma-\beta)$. In fact, in the scale of $\mathcal{FN}_{p,\mu,\infty}^{s}%
$-spaces, the last property occurs only when $s=n-\frac{n-\mu}{p}%
-(2\gamma-\beta)$. So, we obtain that solutions are self-similar, i.e.
invariant under (\ref{scaling}), provided that the symbols $P_{j}$'s are
homogenous of degree $\beta$ and $\theta_{0}$ is homogeneous of degree
$-(2\gamma-\beta)$. Motivated by that, we also analyse the asymptotic
stability of solutions and obtain a basin of attraction around each
self-similar solution.

A space naturally related to $\mathcal{FN}_{p,\mu,\infty}^{s}$ is the
homogeneous Fourier-Besov space $\mathcal{FB}_{p,\infty}^{s}$ which was
introduced by Konieczny-Yoneda \cite{Koni-Yoneda} in order to study the
Navier-Stokes-Coriolis system in $\mathbb{R}^{3}$ with $s=2-\frac{3}{p}$ and
$p>3$. Inspired by \cite{Koni-Yoneda}, we introduce $\mathcal{FN}%
_{p,\mu,\infty}^{s}$ that is larger than the Fourier-Besov space with the same
scaling and seems to be new in analysis of PDEs, at least to the best of our
knowledge. Precisely, we have the continuous inclusion $\mathcal{FB}%
_{p_{2},\infty}^{s_{2}}\subset\mathcal{FN}_{p_{1},\mu,\infty}^{s_{1}}$ for
$1\leq p_{1}\leq p_{2}\leq\infty,$ $0\leq\mu<n$ and $s_{2}+\frac{n}{p_{2}%
}=s_{1}+\frac{n-\mu}{p_{1}}.$ Taking $\mu=0$ in (\ref{FB-norm1}), the Morrey
space $M_{p,\mu}$ coincides with $L^{p},$ and so $\mathcal{FN}_{p,0,\infty
}^{s}=\mathcal{FB}_{p,\infty}^{s}.$ In \cite{Cannone1} the Navier-Stokes
system was studied in $PM^{n-1}$ which is a subspace of $\mathcal{FB}%
_{p,\infty}^{n-1-n/p}$ with norm defined via Fourier transform too. The space
$\mathcal{FN}_{p,\mu,\infty}^{s}$ with $s=n-1-\frac{n-\mu}{p}$ could be
employed to provide a larger initial data class for uniform solvability (with
respect to angular velocity) of the Navier-Stokes-Coriolis system. They also
can be seen as a counterpart in Fourier variables of the homogeneous
Besov-Morrey spaces $\mathcal{N}_{p,\mu,\infty}^{s}$ introduced by
\cite{Kozo1} in order to study the Navier-Stokes equations (see Remark
\ref{Rem-space} for the definition). We also refer the reader to \cite{Mazzu}
for useful properties, bilinear and pseudo-differential estimates in
$\mathcal{N}_{p,\mu,\infty}^{s}$ spaces, and an extension of the analysis to
compact Riemannian manifolds. Due to lack of Hausdorff-Young inequality on
Morrey spaces, it seems that there are no inclusion relations between
$\mathcal{FN}_{p_{1},\mu_{1},\infty}^{s_{1}}$ and $\mathcal{N}_{p_{2},\mu
_{2},\infty}^{s_{2}}$ for $1\leq p_{1},p_{2}<\infty$, $0<\mu_{1},\mu_{2}<n$,
and $s_{2}-\frac{n-\mu_{2}}{p_{2}}=s_{1}+\frac{n-\mu_{1}}{p_{1}}-n$ (same scaling).

Active scalar equations arise in a number of important mathematical toy or
physical models. For instance, in the case $n=1$ we have Burguer's equation
($\beta=1$ and $u=\theta$) and the transport equation
\begin{equation}
\theta_{t}+\left(  \mathcal{H}(\theta)\theta\right)  _{x}+\kappa
(-\Delta)^{\gamma}\theta=0 \label{eqprinc}%
\end{equation}
where $\beta=1$ and $u=\mathcal{H}(\theta)$ is the Hilbert Transform, a
zero-order multiplier operator. For the former example, we refer the reader to
\cite{Kiselev3} and \cite{Dong-2} for results on blow-up, global existence and
regularity of solutions belonging to Lebesgue or smooth spaces. Results on
global existence, self-similarity, finite-time singularity and asymptotic
behavior of solutions for (\ref{eqprinc}) have also been obtained by several
authors, see e.g. \cite{Car-Fer-Pre},\cite{Castro-Cordoba1}%
,\cite{Castro-Cordoba2},\cite{Morlet},\cite{Chae-Cordoba-Fontelos}, and their
references. Replacing in (\ref{eqprinc}) $\mathcal{H}(\theta)\theta$ by
$\theta\Lambda^{\alpha}\mathcal{H}(\theta)$, it leads to%
\begin{equation}
\theta_{t}+\left(  \theta\Lambda^{\alpha}\mathcal{H}(\theta)\right)
_{x}+\kappa(-\Delta)^{\gamma}\theta=0 \label{eqprinc2}%
\end{equation}
which has the form (\ref{dase})-(\ref{velocity}) with $\beta=\alpha+1$ and
$u=\Lambda^{\alpha}\mathcal{H}(\theta)$. Equations (\ref{eqprinc}) and
(\ref{eqprinc2}) are related to models on dislocation dynamics in crystals
with positive $\theta$ representing the number density of fractures per unit
length in the material (see e.g. \cite{Head3,Deslippe}). Self-similar
asymptotic behavior of solutions for (\ref{eqprinc2}) with $u_{0}\in
BUC(\mathbb{R})$ and $\kappa=0$ was proved in \cite{Biler-Karch} showing an
important role of self-similar solutions in the description of asymptotics for
this and related models (see also \cite{Head3}).

In the 2D case, we have the vorticity equation $u=\nabla^{\perp}(-\Delta
)^{-1}\theta$ and SQG equation $u=\nabla^{\perp}((\Delta)^{-\frac{1}{2}}%
\theta)$. The former is a very famous fluid dynamical model while the latter
has been the object of a lot of papers concerning existence, uniqueness,
regularity and asymptotic behavior of solutions in the inviscid case
$\kappa=0$ or in the subcritical ($1/2<\gamma<1$), critical ($\gamma=1/2$) and
supercritical ($\gamma\in(0,1/2)$) ranges (see e.g. \cite{Caff-Vass}%
,\cite{Car-Fer1},\cite{Const2},\cite{Const3},\cite{Const4},\cite{Cordoba1}%
,\cite{Dong-1},\cite{Ju},\cite{Kiselev1},\cite{Kiselev2},\cite{Kiselev3}%
,\cite{NS},\cite{SS}, and their references).

In the 3D case, an example of active scalar PDE comes from magneto-geostrophic
dynamics in a physical situation of fast rotating electrically conducting
fluids, reading as
\begin{equation}
\frac{\partial\theta}{\partial t}+u\nabla_{x}\theta+\left(  -\Delta\right)
^{\gamma}\theta=0,\text{ with }div(u)=0\text{ and }u_{j}=\Sigma_{k=1}%
^{3}\partial_{k}T_{kj}\theta,\label{MG}%
\end{equation}
where $(T_{kj})_{3\times3}$ is a matrix of Calder\'{o}n-Zygmund singular
integral operator (see \cite{Moffat},\cite{Frid-Vicol-3}). Notice that
(\ref{dase})-(\ref{velocity}) recovers (\ref{MG}) with $n=3$ and $\beta=2$. In
fact, (\ref{MG}) is a generalization of the physical case $\gamma=1$ (which is
critical) to include fractional dissipation. An important feature is that the
presence of the underlying magnetic field generates a non-isotropic structure
for the symbol of $T_{kj}$, and then the associated symbols $P_{j}$'s in
(\ref{uks}) are non-radially symmetric (see \cite{Frid-Vicol-2} for an
explicit expression of $T_{kj}$). Global existence of smooth solutions with
$\gamma=1$ and $L^{2}$ initial data can be found in \ \cite{Frid-Vicol-2} by
using De Giorgi techniques in order to reach H\"{o}lder continuity of weak
solutions. See also \cite{Frid-Vicol-1},\cite{Frid-Vicol-4} and their
references for further results on well-posedness, regularity and stability in
a framework of $L^{2}$ or smooth functions.

The dynamical of solutions naturally depend of the coupling between the
velocity and active scalar. So the PDE (\ref{dase}) allows to analyse many
kinds of dynamics by varying the symbol $P_{j}$. The results of \cite{Chae2}
considered $P[\cdot]$ in (\ref{velocity}) such that $P_{j}\in C^{\infty
}(\mathbb{R}^{n}\backslash\{0\}),$ $P_{j}$ is radially symmetric,
nondecreasing in $\left\vert \xi\right\vert ,$ and satisfying a
H\"{o}rmander-Mikhlin type condition. There the authors showed existence of
global solutions in $L^{\infty}((0,\infty);Y),$ where $Y=L^{1}\cap L^{\infty
}\cap B_{q,\infty}^{s,M}$ with $s>1$ and $2\leq q\leq\infty,$ by means of a
priori estimates in Besov spaces and a successive approximation scheme. The
space $B_{q,\infty}^{s,M}$ is an extension of the classical Besov space
$B_{q,\infty}^{s}$ whose the norm increases depending on the growth of $M$. A
technical growth condition depending on $M$ is also assumed for $P_{j}$. Their
results can be applied for several active scalars. For instance, for the
generalized SQG
\begin{equation}
u=\nabla^{\perp}(\Lambda^{\beta-2}\theta)=\Lambda^{\beta-1}(-\mathcal{R}%
_{2}\theta,\mathcal{R}_{1}\theta),\text{ } \label{field11}%
\end{equation}
with $0\leq\beta<2\gamma<1$ and $n=2$. By varying the parameter $\beta$ from
$0$ to $1$, this equation interpolates 2D vorticity and SQG equations, and it
is called modified SQG in the critical case $\beta=2\gamma$. The equation
(\ref{dase}) with (\ref{field11}) has been studied for instance in
\cite{Chae1},\cite{Const1},\cite{Kiselev3},\cite{May1},\cite{Miao1}%
,\cite{Miao-2} where one can find existence and regularity results with data
in Sobolev spaces $H^{m}$ for $m\geq0$. The conditions $\kappa>0$, $\beta
\in\lbrack0,1]$ and $\beta=2\gamma$ were assumed in \cite{Const1}%
,\cite{Kiselev3},\cite{May1},\cite{Miao1}, and $\kappa>0$ and $1\leq
\beta<2\gamma<2$ in \cite{Miao-2}. Also, local well-posedness of $H^{m}%
$-solutions for (\ref{dase})-(\ref{field11}) was proved in the inviscid case
$\kappa=0$ for $\beta\in\lbrack1,2]$ and $m\geq4.$ Another application is for
log-type couplings such as
\begin{align}
P_{j}(\xi)  &  =\left\vert \xi\right\vert ^{\alpha}(\log(1+\left\vert
\xi\right\vert ^{2}))^{\chi},\text{ }\chi>0,\text{ }\label{log-field-2}\\
P_{j}(\xi)  &  =\left\vert \xi\right\vert ^{\alpha}(\log(1+\log(1+\left\vert
\xi\right\vert ^{2})))^{\chi},\chi>0, \label{log-field-3}%
\end{align}
which are interesting in view of the numerical evidences presented in
\cite{Okitani} that, even for $\kappa=0,$ (\ref{dase}) with velocity
\begin{equation}
u=\nabla^{\perp}(\log(I-\Delta))^{\chi}\theta,\text{ }\chi>0,
\label{log-field-1}%
\end{equation}
may be globally well-posed. As a first step for this conjecture, the paper
\cite{Chae1} proved local well-posedness of $H^{4}$-solutions for
(\ref{dase})-(\ref{log-field-1}) with $\gamma>0$. In the case $\alpha=0$ and
$n=2,$ (\ref{log-field-2}) and (\ref{log-field-3}) correspond to \textit{log}
and \textit{log-log} Navier-Stokes, respectively, which are intermediate
models between 2D vorticity and SQG equations. Considering the inviscid case
$\kappa=0$, that is, \textit{log} and \textit{log-log} Euler equations, the
authors of \cite{Chae3} proved global existence results for $0\leq\chi\leq1$
and initial data $\theta_{0}\in$ $L^{1}\cap L^{\infty}\cap B_{q,\infty}^{s}$,
where $s>1$, $q>2$ and $B_{q,\infty}^{s}$ is the inhomogeneous Besov space.

Some symbols with log-type growth also provide examples of slightly
supercritical active scalar equations that behave nicely from the viewpoint of
global existence of smooth solutions. Relying on the method of \cite{Kiselev1}
based on modulus of continuity, the work \cite{Da-Kiselev-Vicol} showed
existence and uniqueness of global smooth solutions for (\ref{dase}) with
$n=2,$ $\gamma=1/2,$ and a multiplier coupling operator
\begin{equation}
u=\nabla^{\perp}\Lambda^{-1}m(\Lambda)\theta, \label{m-coupling}%
\end{equation}
where $m$ is a smooth, radial, nondecreasing function such that $m(\xi)\geq1,$
for all $\xi\in$ $\mathbb{R}^{2}$. The authors also assumed a
H\"{o}rmander-Mikhlin type condition and the growth hypothesis
\begin{equation}
\frac{m(\xi)}{\log\log\left\vert \xi\right\vert }\rightarrow0\text{ as
}\left\vert \xi\right\vert \rightarrow\infty. \label{cond-Vicol}%
\end{equation}
In a certain sense, the coupling (\ref{m-coupling}) is critical at origin and
supercritical at infinity at most by a\textit{ log-log} factor.

In \cite{Fer-Lima}, the authors tackled (\ref{dase}) with high-order couplings
and showed global existence and decay of solutions in the critical Lebesgue
space $L^{\frac{n}{2\gamma-\beta}}(\mathbb{R}^{n})$. There they worked within
a sub-range of $\ 1\leq\beta<2\gamma<n+1$ and considered $P_{j}\in$
$C^{[\frac{n}{2}]+1}(\mathbb{R}^{n}\backslash\{0\})$ such that%

\begin{equation}
\left\vert \frac{\partial^{\alpha}P_{i}}{\partial\xi^{\alpha}}(\xi)\right\vert
\leq C|\xi|^{\beta-|\alpha|}, \label{Pi-cond}%
\end{equation}
for all $\alpha\in(\mathbb{N}\cup\{0\})^{n}$, $|\alpha|\leq\lbrack\frac{n}%
{2}]+1$ and $\xi\neq0$. The condition (\ref{Pi-cond}) is also of
H\"{o}rmander-Mikhlin type and allows non-radial symbols. The solutions in
\cite{Fer-Lima} are $C^{\infty}$-smooth for $t>0,$ but they do not necessarily
belong to $L^{2}$ because $2<\frac{n}{2\gamma-\beta}$ for either $n=2,3$ or
$1/2<\gamma<1.$ The approach in \cite{Fer-Lima} relies on a combination of
time-weighted Kato type norms, scaling arguments, $L^{q}$-maximum principles,
and arguments of the type parabolic De Giorgi-Nash-Moser.

Due to the large number of works cited above, in what follows, we highlight
the novelties of the present paper for the reader convenience. We provide
self-similar solutions and asymptotic behavior results in a new framework
outside $L^{2}$ by means of an approach based on Fourier transform. These
results can be applied in a unified way for the above mentioned couplings in
the subcritical range $\beta<2\gamma$. Our conditions cover symbols
$P_{j}\notin C(\mathbb{R}^{n}\backslash\{0\})$ and non-radial ones like that
of the magneto-geostrophic equation (\ref{MG}). As already pointed out above,
some initial data with large $L^{2}$-norm can be considered because the
smallness condition is with respect to the weak norm of $\mathcal{FN}%
_{p,\mu,\infty}^{s}.$

The outline of this paper is as follows. In section 2, we introduce the
functional setting and recall an useful abstract point fixed lemma. Our
results are stated in section 3. Estimates for the nonlinear term associated
to the mild formulation of (\ref{dase})-(\ref{velocity}) are obtained in
section 4. Finally, the results are proved in section 5.

\section{Preliminares}

\hspace{0.5cm} In this section we introduce Fourier-Besov-Morrey spaces and
recall an abstract fixed point lemma that will be useful for our purposes.

\subsection{Fourier-Besov-Morrey spaces}

\hspace{0.5cm} We start by recalling Morrey spaces $M_{p,\lambda}%
=M_{p,\lambda}(\mathbb{R}^{n})$ (see \cite{Taylor}, \cite{Kato} for more
details). For $1\leq p\leq\infty$ and $0\leq\mu<n,$ the space $M_{p,\mu}$ is
defined as
\begin{equation}
M_{p,\mu}=\left\{  f\in L_{loc}^{p}(\mathbb{R}^{n}):\left\Vert f\right\Vert
_{M_{p,\mu}}<\infty\right\}  , \label{def11}%
\end{equation}
where
\begin{equation}
\left\Vert f\right\Vert _{M_{p,\mu}}=\sup_{x_{0}\in\mathbb{R}^{n}%
,\ 0<R<\infty}\left(  R^{-\frac{\mu}{p}}\left\Vert f\right\Vert _{L^{p}%
(B_{R}(x_{0}))}\right)  \label{def}%
\end{equation}
and $B_{R}(x_{0})\subset$ $\mathbb{R}^{n}$ is the open ball with center
$x_{0}$ and radius $R.$ The space $M_{p,\mu}$ endowed with $\left\Vert
\cdot\right\Vert _{M_{p,\mu}}$ is a Banach space. For $p=\infty,$ $M_{p,\mu}$
becomes $L^{\infty}$. In the case $p=1,$ $M_{p,\mu}$ should be understood as a
space of Radon measures and the $L^{1}$-norm in (\ref{def}) as the total
variation of the measure $f$ on $B_{R}(x_{0}).$

If $1\leq p_{i}\leq\infty$ and $0\leq\mu_{i}<n$ with $\dfrac{1}{p_{3}}%
=\dfrac{1}{p_{1}}+\dfrac{1}{p_{2}}$ and $\dfrac{\mu_{3}}{p_{3}}=\dfrac{\mu
_{1}}{p_{1}}+\dfrac{\mu_{2}}{p_{2}},$ then we have the H\"{o}lder type
inequality
\begin{equation}
\left\Vert fg\right\Vert _{M_{p_{3},\mu_{3}}}\leq\left\Vert f\right\Vert
_{M_{p_{1},\mu_{1}}}\left\Vert g\right\Vert _{M_{p_{2},\mu_{2}}}.
\label{holder}%
\end{equation}
Also, for $1\leq p\leq\infty$ and $0\leq\mu<n,$
\begin{equation}
\left\Vert \varphi\ast g\right\Vert _{M_{p,\mu}}\leq\left\Vert \varphi
\right\Vert _{1}\left\Vert g\right\Vert _{M_{p,\mu}}, \label{young}%
\end{equation}
for all $\varphi\in L^{1}$ and $g\in M_{p,\mu}$

Let us now recall the Littlewood-Paley decomposition (see e.g. \cite{Cannone2}%
,\cite{Kozo1},\cite{Lemarie},\cite{Mazzu}). Let $\varphi\in\mathcal{S}%
(\mathbb{R}^{n})$ be a radially symmetric function with support in $\{\xi
\in\mathbb{R}^{n}:\frac{3}{4}\leq|\xi|\leq\ \frac{8}{3}\}$ and such that
\[
\sum_{k\in\mathbb{Z}}\varphi(2^{-k}\xi)=1,\text{ for all }\xi\neq0.
\]
Consider the family of functions $\{\varphi_{k}\}_{k\in\mathbb{Z}}$ defined by
$\varphi_{k}(\xi)=\varphi(2^{-k}\xi)$, for every $k\in\mathbb{Z}.$ Let
$\mathcal{P}$ be the set of all polynomials with $n$ variables. For
$s\in\mathbb{R}$, $1\leq p\leq\infty$ and $1\leq q\leq\infty$, the homogeneous
Fourier-Besov-Morrey space $\mathcal{FN}_{p,\mu,q}^{s}(\mathbb{R}^{n})$ is the
Banach space of all $f\in\mathcal{S}^{\prime}(\mathbb{R}^{n})/\mathcal{P}$
such that the norm $\left\Vert f\right\Vert _{\mathcal{FN}_{p,\mu,q}^{s}}$ is
finite, where%
\begin{equation}
\left\Vert f\right\Vert _{\mathcal{FN}_{p,\mu,q}^{s}}=%
\begin{cases}
(\sum_{k\in\mathbb{Z}}2^{ksq}\left\Vert \varphi_{k}\widehat{f}\right\Vert
_{M_{p,\mu}}^{q})^{1/q}{\Large ,} & \,\,q<\infty,\\[3mm]%
\sup_{k\in\mathbb{Z}}2^{ks}\left\Vert \varphi_{k}\widehat{f}\right\Vert
_{M_{p,\mu}}, & q=\infty.
\end{cases}
\label{FB-norm1}%
\end{equation}

\begin{remark}
\label{Rem-space}Let us observe that Besov-Morrey spaces $\mathcal{N}%
_{p,\mu,q}^{s}$ studied in \cite{Kozo1},\cite{Mazzu} are obtained by replacing
in (\ref{FB-norm1}) $\varphi_{k}\widehat{f}$ by $(\varphi_{k}\widehat
{f})^{\vee}=[(\varphi_{k})^{\vee}\ast f]$.
\end{remark}

In order to study how the product acts on Fourier-Besov-Morrey space, we need
to recall the Bony's paraproduct formula (see \cite{Bony},\cite{Lemarie}%
,\cite{Mazzu}). Firstly define the localization operators
\begin{equation}
\Delta_{j}f=\varphi_{j}(D)f\text{ and }S_{j}f=\sum_{k\leq j-1}\Delta
_{k}f,\text{ for every }j\in\mathbb{Z}. \label{FB-2}%
\end{equation}
A direct computation gives the equalities
\begin{align}
\Delta_{j}\Delta_{k}f  &  =0,\text{ if }|j-k|\geq2,\label{operator1}\\
\Delta_{j}(S_{k-1}f\Delta_{k}g)  &  =0,\text{ if }|j-k|\geq5.
\label{operator2}%
\end{align}

Given $f\in\mathcal{S}^{\prime}(\mathbb{R}^{n}),$ the paraproduct operator
$T_{f}$ is defined as
\begin{equation}
T_{f}g=\sum_{k\in\mathbb{Z}}S_{k-1}f\Delta_{k}g,\text{ for all }%
g\in\mathcal{S}^{\prime}(\mathbb{R}^{n}). \label{para1}%
\end{equation}
Using (\ref{para1}) we can write
\begin{equation}
fg=T_{f}g+T_{g}f+R(f,g) \label{Bonydec}%
\end{equation}
where
\begin{equation}
R(f,g)=\sum_{k\in\mathbb{Z}}\Delta_{k}f\tilde{\Delta}_{k}g\text{ and }%
\tilde{\Delta}_{k}g=\sum_{|k^{\prime}-k|\leq1}\Delta_{k^{\prime}}g.
\label{para2}%
\end{equation}
It follows from (\ref{operator1})-(\ref{operator2}) that
\begin{equation}
\Delta_{j}(fg)=\sum_{|k-j|\leq4}\Delta_{j}(S_{k-1}f\Delta_{k}g)+\sum
_{|k-j|\leq4}\Delta_{j}(S_{k-1}g\Delta_{k}f)+\sum_{k\geq j-2}\Delta_{j}%
(\Delta_{k}f\widetilde{\Delta}_{k}g). \label{para3}%
\end{equation}

We finish this subsection with a Bernstein type lemma in Fourier variables in
Morrey spaces.

\begin{lemma}
\label{bernstein} Let $1\leq q\leq p\leq\infty,$ $0\leq\mu_{1},\mu_{2}<n,$
$\frac{n-\mu_{2}}{p}\leq\frac{n-\mu_{1}}{q},$ and let $\alpha$ be a
multiindex. If $supp(\hat{f})\subset\{|\xi|\leq A2^{j}\}$ then there is a
constant $C>0$ independent of $f$ and $j$ such that
\begin{equation}
\Vert(i\xi)^{\alpha}\hat{f}\Vert_{M_{q,\mu_{2}}}\leq C2^{j|\alpha
|+j(\frac{n-\mu_{2}}{q}-\frac{n-\mu_{1}}{p})}\Vert\hat{f}\Vert_{M_{p,\mu_{1}}%
}. \label{Bern1}%
\end{equation}

\end{lemma}

\textbf{Proof. \ }\bigskip From the condition on support of $\hat{f}$ and
(\ref{holder}), we obtain
\[
\Vert(i\xi)^{\alpha}\hat{f}\Vert_{M_{q,\mu_{2}}}\leq C2^{j|\alpha|}\Vert
1\cdot\hat{f}\Vert_{M_{q,\mu_{2}}}\leq C2^{j|\alpha|}2^{j(\frac{n-\mu_{2}}%
{q}-\frac{n-\mu_{1}}{p})}\Vert\hat{f}\Vert_{M_{p,\mu_{1}}},
\]
which gives (\ref{Bern1}). \qed

\subsection{\bigskip A fixed point lemma}

\hspace{0.5cm} In order to avoid extensive point fixed arguments we are going
to use the following abstract lemma.

\begin{lemma}
\label{genlem} (see \cite{Lemarie}) Let $X$ be a Banach space with norm
$\Vert\cdot\Vert_{X}$, and $\mathcal{B}:X\times X\rightarrow X$ be a
continuous bilinear map, i.e., there exists $K>0$ such that
\[
\Vert\mathcal{B}(x_{1},x_{2})\Vert_{X}\leq K\text{ }\Vert x_{1}\Vert_{X}\text{
}\Vert x_{2}\Vert_{X},
\]
for all $x_{1},x_{2}\in X$. Given $0<\varepsilon<\frac{1}{4K}$ and $y\in X$
such that $\Vert y\Vert_{X}\leq\varepsilon$, there exists a unique solution
$x\in X$ for the equation $x=y+\mathcal{B}(x,x)$ in the closed ball~$\left\{
x\in X:\left\Vert x\right\Vert _{X}\leq2\varepsilon\right\}  .$ Moreover, the
solution $x$ depends continuously on $y$ in the following sense: If
$\Vert\tilde{y}\Vert_{X}\leq\varepsilon$, $\tilde{x}=\tilde{y}+\mathcal{B}%
(\tilde{x},\tilde{x})$, and $\Vert\tilde{x}\Vert_{X}\leq2\varepsilon$, then
\begin{equation}
\Vert x-\tilde{x}\Vert_{X}\leq\frac{1}{1-4K\varepsilon}\Vert y-\tilde{y}%
\Vert_{X}. \label{lip01}%
\end{equation}

\end{lemma}

\section{Results}

\hspace{0.5cm} In this work we consider solutions for (\ref{dase}%
)-(\ref{velocity}) which satisfy an integral equation come from Duhamel's
principle. This equation reads as
\begin{equation}
\theta(t)=G_{\gamma}(t)\theta_{0}+\mathcal{B}(\theta,\theta)(t), \label{mild}%
\end{equation}
where $G_{\gamma}(t)\theta_{0}=g_{\gamma}(\cdot,t)\ast\theta_{0}$ with
$\hat{g}_{\gamma}(\xi,t)=e^{-t|\xi|^{2\gamma}}$ and
\begin{equation}
\mathcal{B}(\theta,\varphi)(t)=-\int_{0}^{t}G_{\gamma}(t-\tau)(\nabla_{x}%
\cdot(P[\theta]\varphi))(\tau)d\tau. \label{termo bilinear}%
\end{equation}

In what follows, we state our well-posedness and self-similarity results.

\begin{theorem}
\label{teoglobal} Let $n\geq1,$ $\gamma>1/2,$ $0\leq\beta<2\gamma
<\frac{n+\beta+1}{2},$ $\frac{n-\mu}{n+\beta+1-4\gamma}<p\leq\infty,$ and
$0\leq\mu<n.$ Suppose $\theta_{0}\in\mathcal{FN}_{p,\mu,\infty}^{s}$ with
$s=n-\frac{n-\mu}{p}-(2\gamma-\beta).$ Let $K$ be as in Lemma \ref{LemaBil}
and $0<\varepsilon<\frac{1}{4K}$.

\begin{description}
\item[(i)] (Well-posedness) If $\Vert\theta_{0}\Vert_{\mathcal{FN}%
_{p,\mu,\infty}^{s}}\leq\varepsilon$ then (\ref{mild}) has a unique global
solution%
\[
\theta\in BC((0,\infty);\mathcal{FN}_{p,\mu,\infty}^{s}(\mathbb{R}^{n}))
\]
such that $\sup_{t>0}\Vert\theta(t)\Vert_{\mathcal{FN}_{p,\mu,\infty}^{s}}%
\leq2\varepsilon$. The data-solution map $\theta_{0}\rightarrow\theta(x,t)$
from \newline%
\[
\left\{  f\in\mathcal{FN}_{p,\mu,\infty}^{s}:\left\Vert f\right\Vert
_{\mathcal{FN}_{p,\mu,\infty}^{s}}\leq\varepsilon\right\}  \text{ to
}BC((0,\infty);\mathcal{FN}_{p,\mu,\infty}^{s}(\mathbb{R}^{n}))
\]
is Lipschitz continuous. Moreover $\theta(x,t)\rightharpoonup\theta_{0}$ in
$\mathcal{S}^{\prime}(\mathbb{R}^{n})$ as $t\rightarrow0^{+}$.

\item[(ii)] (Self-Similarity) If the symbol $P_{j}(\xi)$ is homogeneous of
degree $\beta$, for $j=1,2,...,n,$ and if $\theta_{0}$ is a homogeneous
distribution of degree $-(2\gamma-\beta)$ then the solution $\theta$ is
self-similar, that is,
\[
\theta=\theta_{\lambda}:=\lambda^{2\gamma-\beta}\theta(\lambda x,\lambda
^{2\gamma}t),\text{ for all }\lambda>0\text{.}%
\]

\end{description}
\end{theorem}

\bigskip

\begin{remark}
(Local-in-time solutions) Assume the same conditions on $n,\beta,\gamma,p,\mu$
given in Theorem \ref{teoglobal} and let $s>n-\frac{n-\mu}{p}-(2\gamma-\beta)$
and $\theta_{0}\in\mathcal{FN}_{p,\mu,\infty}^{s}.$ A local-in-time
well-posedness result in $BC((0,T);\mathcal{FN}_{p,\mu,\infty}^{s}%
(\mathbb{R}^{n}))$ for (\ref{mild}) could be proved by assuming a small
condition on $T>0$ and regardless the size of the initial data. Again, the
solution $\theta$ would satisfy the initial condition in the sense of
distributions, that is, $\theta(x,t)\rightharpoonup\theta_{0}$ in
$\mathcal{S}^{\prime}(\mathbb{R}^{n})$ as $t\rightarrow0^{+}$.
\end{remark}

\bigskip

\begin{remark}
(Data with large $L^{2}$-norm) \label{Rem2} Let $\theta_{0},\psi_{0}%
\in\mathcal{S}^{\prime}(\mathbb{R}^{n})$ be such that%
\begin{equation}
\hat{\theta}_{0}(\xi)=\delta\left\vert \xi\right\vert ^{-(n-(2\gamma-\beta
))}1_{\left\vert \xi\right\vert <R_{1}}\text{ and \ }\hat{\psi}_{0}%
(\xi)=\delta\left\vert \xi\right\vert ^{-(n-(2\gamma-\beta))}1_{\left\vert
\xi\right\vert >R_{2}}, \label{data-Rem}%
\end{equation}
where $\delta>0$ and $1_{A}$ stands for the characteristic function of set
$A$. Recalling that $(\left\vert \xi\right\vert ^{-(n-(2\gamma-\beta))}%
)^{\vee}=C\left\vert x\right\vert ^{-(2\gamma-\beta)}$, it follows from
(\ref{FB-norm1}) that
\[
\left\Vert \theta_{0}\right\Vert _{\mathcal{FN}_{p,\mu,\infty}^{s}},\left\Vert
\psi_{0}\right\Vert _{\mathcal{FN}_{p,\mu,\infty}^{s}}\leq\left\Vert \delta
C\left\vert x\right\vert ^{-(2\gamma-\beta)}\right\Vert _{\mathcal{FN}%
_{p,\mu,\infty}^{s}}=C\delta,
\]
where $C>0$ is independent of $R_{1},R_{2}$. We can apply Theorem
\ref{teoglobal} with $\theta_{0}$ and $\psi_{0}$ for all $R_{1},R_{2}>0$ by
choosing $\delta>0$ in such a way that $\delta=\frac{\varepsilon}{C}.$ On the
other hand, if $n+2\beta<4\gamma,$ then $\theta_{0}\in L^{2}$ with $\left\Vert
\theta_{0}\right\Vert _{L^{2}}^{2}=C(\delta)R_{1}^{4\gamma-(n+2\beta)}$ for
all $R_{1}>0,$ and if $n+2\beta>4\gamma$ then $\psi_{0}\in L^{2}$ with
$\left\Vert \psi_{0}\right\Vert _{L^{2}}^{2}=C(\delta)R_{2}^{4\gamma
-(n+2\beta)}$, for all $R_{2}>0$. Then $\theta_{0}$ and $\psi_{0}$ can have
arbitrarily large $L^{2}$-norm by making $R_{1}\rightarrow\infty$ and
$R_{2}\rightarrow0,$ respectively.
\end{remark}

Now we present a result of asymptotic stability of solutions in the framework
of Fourier-Besov-Morrey spaces. Here we follow ideas from \cite{Cannone1}
where the authors studied Navier-Stokes equations in the framework of $PM^{a}$-spaces.

\begin{theorem}
(Asymptotic stability)\label{asymp} Under the hypotheses of Theorem
\ref{teoglobal}. Assume that $\theta$ and $\phi$ are solutions for
(\ref{mild}) given by Theorem \ref{teoglobal} with small initial data
$\theta_{0}$ and $\phi_{0}\in\mathcal{FN}_{p,\mu,\infty}^{s}(\mathbb{R}^{n})$,
respectively. We have that
\begin{equation}
\lim_{t\rightarrow\infty}\Vert\theta(\cdot,t)-\phi(\cdot,t)\Vert
_{\mathcal{FN}_{p,\mu,\infty}^{s}}=0 \label{Asymp1}%
\end{equation}
if and only if%
\begin{equation}
\lim_{t\rightarrow\infty}\Vert G_{\gamma}(t)(\theta_{0}-\phi_{0}%
)\Vert_{\mathcal{FN}_{p,\mu,\infty}^{s}}=0. \label{cond1}%
\end{equation}

\end{theorem}

\begin{remark}
Theorem \ref{asymp} provides a class of solutions asymptotically self-similar
at infinity. In fact, when $\phi_{0}=$ $\theta_{0}+\varphi$ with $\varphi
\in\mathcal{S}(\mathbb{R}^{n})$ and $\theta_{0}$ is homogeneous of degree
$-(2\gamma-\beta),$ {the corresponding solution }$\phi(x,t)$ is attracted to
the self-similar solution $\theta(x,t)$ in the sense of (\ref{Asymp1}).
\end{remark}

\section{Bilinear estimate}

\hspace{0.5cm} In order to perform a fixed point argument in
$BC((0,T);\mathcal{FN}_{p,\mu,\infty}^{s})$, we need to obtain estimates for
the bilinear part $\mathcal{B}(\cdot,\cdot)$ of (\ref{mild}) on this space.

\begin{lemma}
\label{LemaBil} Under the hypothesis of Theorem \ref{teoglobal}, there exists
a constant $K>0$ such that
\begin{equation}
\sup_{t>0}\Vert\mathcal{B}(\theta,\phi)\Vert_{\mathcal{FN}_{p,\mu,\infty}^{s}%
}\leq K\sup_{t>0}\Vert\theta\Vert_{\mathcal{FN}_{p,\mu,\infty}^{s}}\sup
_{t>0}\Vert\phi\Vert_{\mathcal{FN}_{p,\mu,\infty}^{s}},\label{bilinear}%
\end{equation}
for all $\theta,\,\phi\in L^{\infty}((0,\infty);\mathcal{FN}_{p,\mu,\infty
}^{s}(\mathbb{R}^{n})).$
\end{lemma}

\textbf{Proof. \ }\bigskip Let us first prove that there exists $K_{1}>0$ such
that
\begin{equation}
\sup_{t>0}\Vert\mathcal{B}(\theta,\phi)\Vert_{\mathcal{FN}_{p,\mu,\infty}^{s}%
}\leq K_{1}\sup_{t>0}\Vert P[\theta]\phi\Vert_{\mathcal{FN}_{p,\mu,\infty
}^{s-2\gamma+1}} \label{aux-bili1}%
\end{equation}
For this let $k\in\mathbb{Z}$ be fixed. Since $supp(\varphi_{k})\subset
\{\xi\in\mathbb{R}^{n}:2^{k-1}\leq|\xi|\leq2^{k+1}\}$ we have
\begin{align}
2^{ks}\Vert\varphi_{k}\widehat{\mathcal{B}(\theta,\phi)}\Vert_{M_{p,\mu}}  &
\leq2^{ks}\int_{0}^{t}\Vert\,|\xi|\,e^{-(t-\tau)\left\vert \xi\right\vert
^{2\gamma}}(\widehat{P[\theta]\phi})(\xi,\tau)\varphi_{k}(\xi)\Vert_{M_{p,\mu
}}\text{ }d\tau\nonumber\\
&  \leq\int_{0}^{t}2^{k+1}e^{-(t-\tau)2^{2\gamma(k-1)}}2^{(2\gamma
-1)k}2^{(s-2\gamma+1)k}\Vert(\widehat{P[\theta]\phi})(\cdot,\tau)\varphi
_{k}(\cdot)\Vert_{M_{p,\mu}}d\tau\nonumber\\
&  \leq2^{2\gamma+1}\int_{0}^{t}2^{2\gamma(k-1)}e^{-(t-\tau)2^{2\gamma(k-1)}%
}2^{(s-2\gamma+1)k}\Vert(\widehat{P[\theta]\phi})(\cdot,\tau)\varphi_{k}%
(\cdot)\Vert_{M_{p,\mu}}\text{ }d\tau\nonumber\\
&  \leq C\int_{0}^{t}2^{2\gamma(k-1)}e^{-(t-\tau)2^{2\gamma(k-1)}}%
\Vert(P[\theta]\phi)(\cdot,\tau)\Vert_{\mathcal{FN}_{p,\mu,\infty}%
^{s-2\gamma+1}}d\tau\label{aux-bili1-2}\\
&  \leq K_{1}\sup_{t>0}\Vert(P[\theta]\phi)(t)\Vert_{\mathcal{FN}%
_{p,\mu,\infty}^{s-2\gamma+1}}. \label{aux-bili2}%
\end{align}
Taking the supremum over $k\in\mathbb{Z}$ and afterwards over $t>0$, we obtain
(\ref{aux-bili1}).

It remains to prove the product estimate
\begin{equation}
\sup_{t>0}\Vert P[\theta]\phi\Vert_{\mathcal{FN}_{p,\mu,\infty}^{s-2\gamma+1}%
}\leq K_{2}\sup_{t>0}\Vert\theta\Vert_{\mathcal{FN}_{p,\mu,\infty}^{s}}%
\sup_{t>0}\Vert\phi\Vert_{\mathcal{FN}_{p,\mu,\infty}^{s}}.\label{aux-bili4}%
\end{equation}
From estimate (\ref{est velocity}), we get%
\begin{equation}
\Vert\varphi_{k}(\xi)\hat{u}(\xi,t)\Vert_{M_{p,\mu}}\leq C2^{k(\beta-1)}%
\Vert\varphi_{k}(\xi)\hat{\theta}(\xi,t)\Vert_{M_{p,\mu}},\label{aux-velo}%
\end{equation}
because $supp(\varphi_{k})\subseteq\{\xi\in\mathbb{R}^{n};2^{k-1}\leq|\xi
|\leq2^{k+1}\}.$

From definition (\ref{FB-norm1}) and (\ref{FB-2}) we have that%
\[
\Vert P[\theta]\phi\Vert_{\mathcal{FN}_{p,\mu,\infty}^{s-2\gamma+1}}%
=\sup_{j\in\mathbb{Z}}2^{j(s-2\gamma+1)}\Vert\widehat{\Delta_{j}(P[\theta
]\phi)}\Vert_{M_{p,\mu}}.
\]
Let $j\in\mathbb{Z}$ be fixed. Using the decomposition (\ref{Bonydec}) with
$f$ and $g$ replaced by $P[\theta]$ and $\phi$, respectively, and taking the
norm $\Vert\cdot\Vert_{M_{p,\mu}},$ we get
\begin{align}
\Vert\varphi_{j}\widehat{P[\theta]\phi}\Vert_{M_{p,\mu}} &  \leq
\sum_{|k-j|\leq4}\Vert\varphi_{j}(\widehat{S_{k-1}P[\theta]})\ast(\varphi
_{k}\hat{\phi})\Vert_{M_{p,\mu}}+\sum_{|k-j|\leq4}\Vert\varphi_{j}%
(\widehat{S_{k-1}\phi})\ast(\varphi_{k}\widehat{P[\theta]})\Vert_{M_{p,\mu}%
}\nonumber\\
&  +\sum_{k\geq j-2}\Vert\varphi_{j}(\varphi_{k}\widehat{P[\theta]}%
)\ast(\widetilde{\varphi}_{k}\hat{\phi})\Vert_{M_{p,\mu}}=I_{1}+I_{2}%
+I_{3}.\label{aux-bili3}%
\end{align}
In order to estimate $I_{1}$ we use Young's inequality in Morrey spaces
(\ref{young}) to obtain
\begin{equation}
2^{j(s-2\gamma+1)}I_{1}\leq2^{j(s-2\gamma+1)}\sum_{|k-j|\leq4}\Vert
\widehat{S_{k-1}P[\theta]}\Vert_{L^{1}}\Vert\varphi_{k}\hat{\phi}%
\Vert_{M_{p,\mu}}.\label{aux-prod1}%
\end{equation}
Recalling the definition of $S_{j}$, Lemma \ref{bernstein} with $\left\vert
\alpha\right\vert =0$ implies that
\begin{equation}
\Vert\widehat{S_{k-1}P[\theta]}\Vert_{L^{1}}\leq\sum_{k^{\prime}<k}%
\Vert\varphi_{k^{\prime}}\widehat{P[\theta]}\Vert_{L^{1}}\leq C\sum
_{k^{\prime}<k}2^{k^{\prime}(n-\frac{n-\mu}{p})}\Vert\varphi_{k^{\prime}%
}\widehat{P[\theta]}\Vert_{M_{p,\mu}}.\label{aux-prod2}%
\end{equation}
Since $n-\frac{n-\mu}{p}=s+2\gamma-\beta$ and $2^{k}\sim2^{j}$ when
$|k-j|\leq4$, it follows from (\ref{aux-prod2}) that
\begin{align}
&  2^{j(s-2\gamma+1)}\sum_{|k-j|\leq4}\Vert\widehat{S_{k-1}P[\theta]}%
\Vert_{L^{1}}\Vert\varphi_{k}\hat{\phi}\Vert_{M_{p,\mu}}\nonumber\\
&  \leq2^{j(s-2\gamma+1)}\sum_{|k-j|\leq4}\sum_{k^{\prime}<k}2^{k^{\prime
}(s+2\gamma-\beta)}\Vert\varphi_{k^{\prime}}\widehat{P[\theta]}\Vert
_{M_{p,\mu}}\Vert\varphi_{k}\hat{\phi}\Vert_{M_{p,\mu}}\label{aux-prod3}\\
&  \leq2^{j(s-2\gamma+1)}C\sum_{|k-j|\leq4}\Vert\varphi_{k}\widehat{\phi}%
\Vert_{M_{p,\mu}}\sum_{k^{\prime}<k}2^{k^{\prime}(2\gamma-1)}2^{k^{\prime
}(s-\beta+1)}2^{k^{\prime}(\beta-1)}\Vert\varphi_{k^{\prime}}\widehat{\theta
}\Vert_{M_{p,\mu}}\nonumber\\
&  \leq2^{j(s-2\gamma+1)}C\left(  \sum_{|k-j|\leq4}\Vert\varphi_{k}%
\widehat{\phi}\Vert_{M_{p,\mu}}\sum_{k^{\prime}<k}2^{k^{\prime}(2\gamma
-1)}\right)  \left(  \sup_{k^{\prime}\in\mathbb{Z}}2^{k^{\prime}s}\Vert
\varphi_{k^{\prime}}\widehat{\theta}\Vert_{M_{p,\mu}}\right)
\label{aux-prod333}\\
&  \leq C\Vert\theta\Vert_{\mathcal{FN}_{p,\mu,\infty}^{s}}\sum_{|k-j|\leq
4}2^{k(s-2\gamma+1)}\Vert\varphi_{k}\widehat{\phi}\Vert_{M_{p,\mu}%
}2^{k(2\gamma-1)}\nonumber\\
&  \leq C\Vert\theta\Vert_{\mathcal{FN}_{p,\mu,\infty}^{s}}\sum_{|k-j|\leq
4}2^{ks}\Vert\varphi_{k}\widehat{\phi}\Vert_{M_{p,\mu}}\nonumber\\
&  \leq C\sup_{t>0}\Vert\theta\Vert_{\mathcal{FN}_{p,\mu,\infty}^{s}}%
\sup_{t>0}\Vert\phi\Vert_{\mathcal{FN}_{p,\mu,\infty}^{s}},\label{aux-prod4}%
\end{align}
where we have used (\ref{aux-velo}) in (\ref{aux-prod3}), and $\gamma>1/2$ in
(\ref{aux-prod333}). Inserting (\ref{aux-prod4}) into (\ref{aux-prod1}) we
obtain
\begin{equation}
2^{j(s-2\gamma+1)}I_{1}\leq C\sup_{t>0}\Vert\theta\Vert_{\mathcal{FN}%
_{p,\mu,\infty}^{s}}\sup_{t>0}\Vert\phi\Vert_{\mathcal{FN}_{p,\mu,\infty}^{s}%
}.\label{aux-1}%
\end{equation}
By using that $\beta<2\gamma$, similarly one can show that
\begin{equation}
2^{j(s-2\gamma+1)}I_{2}\leq C\sup_{t>0}\Vert\theta\Vert_{\mathcal{FN}%
_{p,\mu,\infty}^{s}}\sup_{t>0}\Vert\phi\Vert_{\mathcal{FN}_{p,\mu,\infty}^{s}%
}.\label{aux-2}%
\end{equation}

Now we deal with the third term in (\ref{aux-bili3}). Note that $s-2\gamma
+1>0$ due to the condition $\frac{n-\mu}{n+\beta+1-4\gamma}<p$. Again
employing Young's inequality and Lemma \ref{bernstein}, we get
\begin{align}
2^{j(s-2\gamma+1)}I_{3} &  \leq2^{j(s-2\gamma+1)}\sum_{k\geq j-2}\Vert
\varphi_{k}\widehat{P[\theta]}\Vert_{L^{1}}\Vert\tilde{\varphi}_{k}%
\widehat{\phi}\Vert_{M_{p,\mu}}\nonumber\\
&  \leq\sum_{k\geq j-2}2^{j(s-2\gamma+1)}2^{k(\beta-1)+k(n-\frac{n-\mu}{p}%
)}\Vert\varphi_{k}\widehat{\theta}\Vert_{M_{p,\mu}}\Vert\tilde{\varphi}%
_{k}\widehat{\phi}\Vert_{M_{p,\mu}}\nonumber\\
&  \leq\sum_{k\geq j-2}2^{j(s-2\gamma+1)}2^{(s+2\gamma-1)k}\Vert\varphi
_{k}\widehat{\theta}\Vert_{M_{p,\mu}}\Vert\tilde{\varphi}_{k}\widehat{\phi
}\Vert_{M_{p,\mu}}\nonumber\\
&  \leq C\sum_{k\geq j-2}2^{(j-k)(s-2\gamma+1)}\left(  2^{ks}\Vert\varphi
_{k}\widehat{\theta}\Vert_{M_{p,\mu}}\right)  \left(  2^{ks}\Vert
\tilde{\varphi}_{k}\hat{\phi}\Vert_{M_{p,\mu}}\right)  \nonumber\\
&  \leq C\sup_{t>0}\Vert\theta\Vert_{\mathcal{FN}_{p,\mu,\infty}^{s}}%
\sup_{t>0}\Vert\phi\Vert_{\mathcal{FN}_{p,\mu,\infty}^{s}},\label{aux-3}%
\end{align}
because $\sum_{k\geq j-2}2^{(j-k)(s-2\gamma+1)}<\infty$. The estimates
(\ref{aux-1}), (\ref{aux-2}) and (\ref{aux-3}) yield (\ref{aux-bili4}), as required.\qed

\section{Proof of Theorems}

\hspace{0.5cm} In this section we prove Theorems \ref{teoglobal} and
\ref{asymp}. The proof of the former is the subject of the next subsection.

\subsection{Proof of Theorem \ref{teoglobal}}

\

\textbf{Part (i):} Let us define the Banach space
\[
\mathcal{X}=BC((0,\infty);\mathcal{FN}_{p,\mu,\infty}^{s}(\mathbb{R}^{n}))
\]
with norm given by
\[
\left\Vert \theta\right\Vert _{\mathcal{X}}=\sup_{t>0}\Vert\theta
(\cdot,t)\Vert_{\mathcal{FN}_{p,\mu,\infty}^{s}}.
\]
Rewriting Lemma \ref{LemaBil} with the norm $\left\Vert \cdot\right\Vert
_{\mathcal{X}}$, we obtain
\begin{equation}
\Vert\mathcal{B}(\theta,\phi)\Vert_{\mathcal{X}}\leq K\Vert\theta
\Vert_{\mathcal{X}}\Vert\phi\Vert_{\mathcal{X}}. \label{proof-1}%
\end{equation}
The linear part in (\ref{mild}) can be handle as
\begin{align}
\left\Vert G_{\gamma}(t)\theta_{0}\right\Vert _{\mathcal{X}}  &  =\sup
_{t>0}\left(  \sup_{k\in\mathbb{Z}}2^{ks}\Vert\varphi_{k}e^{-t\left\vert
\xi\right\vert ^{2\gamma}}\hat{\theta}_{0}\Vert_{M_{p,\mu}}\right) \nonumber\\
&  \leq\sup_{t>0}\left(  \sup_{k\in\mathbb{Z}}2^{ks}\Vert\varphi_{k}%
\hat{\theta}_{0}\Vert_{M_{p,\mu}}\right) \nonumber\\
&  =\Vert\theta_{0}\Vert_{\mathcal{FN}_{p,\mu,\infty}^{s}}. \label{proof-2}%
\end{align}

So, taking $0<\varepsilon<\frac{1}{4K}$ and $\Vert\theta_{0}\Vert
_{\mathcal{FN}_{p,\mu,\infty}^{s}}\leq\varepsilon,$ Lemma \ref{genlem}
together with the estimates (\ref{proof-1}) and (\ref{proof-2}) assures that
there is a unique solution $\theta\in\mathcal{X}$ for (\ref{mild}) such that
$\Vert\theta\Vert_{\mathcal{X}}\leq2\varepsilon$. The Lipschitz continuity is
obtained from (\ref{lip01}). The weak convergence as $t\rightarrow0^{+}$
follows by standard arguments and the reader is referred to \cite{Yamazaki}. \qed

\textbf{Part (ii):} Since the symbol $P_{j}$ is homogeneous of degree $\beta$,
for $j=1,2,...,n,$ and $\theta_{0}$ homogeneous of degree $-(2\gamma-\beta)$,
a simple computation shows that $\theta_{\lambda}=\lambda^{2\gamma-\beta
}\theta(\lambda x,\lambda^{2\gamma}t)$ verifies (\ref{mild}), for all
$\lambda>0$, provided that $\theta$ does so. Moreover, due to the scaling
invariance of the $\mathcal{X}$-norm, it follows that $\Vert\theta_{\lambda
}\Vert_{\mathcal{X}}=\Vert\theta\Vert_{\mathcal{X}}\leq2\varepsilon.$ The
uniqueness result contained in item $(i)$ gives us that $\theta_{\lambda
}=\theta$, for all $\lambda>0.$ \qed

\subsection{Proofs of Theorem \ref{asymp}}

\hspace{0.5cm} We prove that (\ref{cond1}) implies (\ref{Asymp1}). The
converse follows similarly and is left to the reader. Subtracting the
corresponding integral equations verified by $\theta$ and $\phi$ and
afterwards computing the $\mathcal{FN}_{p,\mu,\infty}^{s}$-norm, we obtain
\begin{align}
\Vert\theta(t)-\phi(t)\Vert_{\mathcal{FN}_{p,\mu,\infty}^{s}}  &  \leq\Vert
G_{\gamma}(t)(\theta_{0}-\phi_{0})\Vert_{\mathcal{FN}_{p,\mu,\infty}^{s}%
}+\Vert B(\theta-\phi,\theta)+B(\phi,\theta-\phi)\Vert_{\mathcal{FN}%
_{p,\mu,\infty}^{s}}\nonumber\\
&  \leq\Vert G_{\gamma}(t)(\theta_{0}-\phi_{0})\Vert_{\mathcal{FN}%
_{p,\mu,\infty}^{s}}\nonumber\\
&  +\left\Vert \int_{0}^{\delta t}G_{\gamma}(t-\tau)\nabla_{x}\cdot\left(
P[\theta-\phi]\theta(\tau)+P[\phi](\theta-\phi)(\tau)\right)  d\tau\right\Vert
_{\mathcal{FN}_{p,\mu,\infty}^{s}}\nonumber\\
&  +\left\Vert \int_{\delta t}^{t}G_{\gamma}(t-\tau)\nabla_{x}\cdot\left(
P[\theta-\phi]\theta(\tau)+P[\phi](\theta-\phi)(\tau)\right)  d\tau\right\Vert
_{\mathcal{FN}_{p,\mu,\infty}^{s}}\nonumber\\
&  :=I_{0}+I_{1}+I_{2}, \label{aux-asymp1}%
\end{align}
where $\delta>0$ will be chosen later. Proceeding as in (\ref{aux-bili1-2})
and using (\ref{aux-bili4}), we have that
\begin{align}
I_{1}  &  \leq C\sup_{k\in\mathbb{Z}}\int_{0}^{\delta t}2^{2\gamma
(k-1)}e^{-(t-\tau)2^{2\gamma(k-1)}}\Vert\theta(\tau)-\phi(\tau)\Vert
_{\mathcal{FN}_{p,\mu,\infty}^{s}}\left(  \Vert\theta(\tau)\Vert
_{\mathcal{FN}_{p,\mu,\infty}^{s}}+\Vert\phi(\tau)\Vert_{\mathcal{FN}%
_{p,\mu,\infty}^{s}}\right)  d\tau\nonumber\\
&  \leq4\varepsilon K\sup_{k\in\mathbb{Z}}\int_{0}^{\delta t}2^{2\gamma
(k-1)}e^{-(t-\tau)2^{2\gamma(k-1)}}\Vert\theta(\tau)-\phi(\tau)\Vert
_{\mathcal{FN}_{p,\mu,\infty}^{s}}d\tau, \label{aux-asym2}%
\end{align}
where $K>0$ is as in Lemma \ref{LemaBil}, and above we have used that
\begin{equation}
\sup_{t>0}\Vert\theta(t)\Vert_{\mathcal{FN}_{p,\mu,\infty}^{s}}\leq
2\varepsilon\text{ and }\sup_{t>0}\Vert\phi(t)\Vert_{\mathcal{FN}%
_{p,\mu,\infty}^{s}}\leq2\varepsilon. \label{bola}%
\end{equation}
Making the change of variables $\tau=tz$ in (\ref{aux-asym2}) and using the
equality
\[
\sup_{k\in\mathbb{Z}}2^{2\gamma(k-1)}e^{-t(1-\tau)2^{2\gamma(k-1)}}=\frac
{C}{t(1-\tau)},
\]
it follows that
\begin{align}
I_{1}  &  \leq4\varepsilon K\sup_{k\in\mathbb{Z}}\int_{0}^{\delta t}%
2^{2\gamma(k-1)}e^{-(t-\tau)2^{2\gamma(k-1)}}\Vert\theta(\tau)-\phi(\tau
)\Vert_{\mathcal{FN}_{p,\mu,\infty}^{s}}d\tau\nonumber\\
&  =\sup_{k\in\mathbb{Z}}\int_{0}^{\delta}t2^{2\gamma(k-1)}e^{-t(1-\tau
)2^{2\gamma(k-1)}}\Vert\theta(t\tau)-\phi(t\tau)\Vert_{\mathcal{FN}%
_{p,\mu,\infty}^{s}}d\tau\nonumber\\
&  \leq4\varepsilon KC\int_{0}^{\delta}(1-\tau)^{-1}\Vert\theta(t\tau
)-\phi(t\tau)\Vert_{\mathcal{FN}_{p,\mu,\infty}^{s}}d\tau. \label{aux-asymp2}%
\end{align}
For $I_{2},$ we have that
\begin{align}
I_{2}  &  \leq C\sup_{k\in\mathbb{Z}}\int_{\delta t}^{t}2^{2\gamma
(k-1)}e^{-(t-\tau)2^{2\gamma(k-1)}}\Vert\theta(\tau)-\phi(\tau)\Vert
_{\mathcal{FN}_{p,\mu,\infty}^{s}}\left(  \Vert\theta(\tau)\Vert
_{\mathcal{FN}_{p,\mu,\infty}^{s}}+\Vert\phi(\tau)\Vert_{\mathcal{FN}%
_{p,\mu,\infty}^{s}}\right)  d\tau\nonumber\\
&  \leq4\varepsilon K\left(  \sup_{k\in\mathbb{Z}}\int_{\delta t}%
^{t}2^{2\gamma(k-1)}e^{-(t-\tau)2^{2\gamma(k-1)}}d\tau\right)  \left(
\sup_{\delta t<\tau<t}\Vert\theta(\tau)-\phi(\tau)\Vert_{\mathcal{FN}%
_{p,\mu,\infty}^{s}}\right) \nonumber\\
&  \leq4\varepsilon K\sup_{\delta t<\tau<t}\Vert\theta(\tau)-\phi(\tau
)\Vert_{\mathcal{FN}_{p,\mu,\infty}^{s}}. \label{aux-asymp3}%
\end{align}
Finally, taking
\[
L=\limsup_{t\rightarrow\infty}\Vert\theta(\cdot,t)-\phi(\cdot,t)\Vert
_{\mathcal{FN}_{p,\mu,\infty}^{s}},
\]
notice that $0\leq L\leq4\varepsilon$ in view of (\ref{bola}). Computing the
$\lim\sup$ in (\ref{aux-asymp1}), the inequalities (\ref{aux-asymp2}) and
(\ref{aux-asymp3}) give us
\[
L\leq\left(  C4\varepsilon K\log(\frac{1}{1-\delta})+4\varepsilon K\right)  L
\]
which implies $L=0,$ because we can choose $\delta>0$ sufficiently small so
that
\[
C4\varepsilon K\log(\frac{1}{1-\delta})+4\varepsilon K<1.
\]

\qed

\end{document}